# NON-UNIFORM IN TIME ROBUST GLOBAL ASYMPTOTIC OUTPUT STABILITY FOR DISCRETE-TIME SYSTEMS

**Iasson Karafyllis**
Division of Mathematics, Dept. of Economics, University of Athens
8 Pesmazoglou Str., 10559, Athens, Greece
email: ikarafil@econ.uoa.gr

**Abstract**
In this paper the notions of non-uniform in time Robust Blobal Asymptotic Output Stability (RGAOS) and Input-to-Output Stability (IOS) for discrete-time systems are studied. Characterizations as well as links between these notions are provided. Particularly, it is shown that a discrete-time system with continuous dynamics satisfies the non-uniform in time IOS property if and only if the corresponding unforced system is non-uniformly in time RGAOS. Necessary and sufficient conditions for the solvability of the Robust Output Feedback Stabilization (ROFS) problem are also given.

**Keywords:** Discrete-Time Systems, Lyapunov Functions, Feedback Stabilization.

## 1. Introduction

In this paper we study discrete-time time-varying systems with outputs:

$$x(t+1) = f(t, d(t), x(t))$$
$$Y(t) = H(t, x(t))$$
$$x(t) \in \mathcal{X}, d(t) \in D, Y(t) \in \mathcal{Y}, t \in Z^+$$
(1.1)

where $\mathcal{X}, \mathcal{Y}$ is a pair of normed linear spaces, $D$ is the set of disturbances (or time-varying parameters) and $f : Z^+ \times D \times \mathcal{X} \to \mathcal{X}$, $H : Z^+ \times \mathcal{X} \to \mathcal{Y}$ with $f(t,d,0) = 0$ and $H(t,0) = 0$ for all $(t,d) \in Z^+ \times D$.

The notion of non-uniform in time Robust Global Asymptotic Output Stability (RGAOS) was introduced and studied in [11] for a wide class of systems, including discrete-time time-varying systems (1.1). In this paper we present Lyapunov-like conditions for non-uniform in time RGAOS for discrete-time systems of the form (1.1). Our results are based on the Lyapunov characterization of Robust Global Asymptotic Stability (RGAS) given in [12] for discrete-time systems and are parallel to the results provided in [15], [16] and [6] for uniform global asymptotic stability with respect to closed sets in a finite-dimensional state space.

The notion of non-uniform in time RGAOS is closely related to the notion of non-uniform in time Input-to-Output Stability (IOS) introduced in [11] for a wide class of systems, including discrete-time time-varying systems of the form:

$$x(t+1) = f(t, d(t), x(t), u(t))$$
$$x(t) \in \mathcal{X}, d(t) \in D, u(t) \in U, t \in Z^+$$
(1.2a)

$$Y(t) = H(t, x(t)), Y(t) \in \mathcal{Y}$$
(1.2b)

where $\mathcal{X}, \mathcal{Y}, U$ is a triplet of normed linear spaces, $D$ is the set of disturbances (or time-varying parameters) and $f : Z^+ \times D \times \mathcal{X} \times U \to \mathcal{X}$, $H : Z^+ \times \mathcal{X} \to \mathcal{Y}$ are mappings with $f(t,d,0,0) = 0$ and $H(t,0) = 0$ for all $(t,d) \in Z^+ \times D$.

The notion of non-uniform in time IOS introduced in [11] extends the notion of uniform in time Input-to-State Stability (ISS) for discrete-time systems introduced in [4] and further studied in [5], [7], [8]. It is also an extension of the notion of uniform in time IOS introduced in [24], [25] and [9] for continuous time systems. In this paper we derive characterizations of non-uniform in time IOS. Moreover, it is shown that a discrete-time time-varying system (1.2) with continuous dynamics satisfies the non-uniform in time Input-to-Output Stability property if and only if the "unforced" system (1.2), i.e., system (1.2) with $u(t) \equiv 0$

$$x(t+1) = f(t, d(t), x(t), 0)$$
$$Y(t) = H(t, x(t))$$
$$x(t) \in \mathcal{X}, d(t) \in D, Y(t) \in \mathcal{Y}, t \in Z^+$$
(1.3)

is non-uniformly in time Robustly Globally Asymptotically Output Stable. This result is important, since it shows that continuity of the dynamics guarantees useful robustness properties.

Discrete-time control systems of the form (1.2) arise naturally in applications. For example in [21], the stability of infinite-dimensional discrete-time systems is studied. In the present paper we focus on the Robust Output Feedback Stabilization problem for (1.2), i.e., the problem of the stabilization of the output (1.2b) of the time-varying discrete-time system (1.2a) by means of:

**(i)** a time-varying output feedback law $u(t) = k(t, y(t))$ (static ROFS problem)

**(ii)** a dynamic time-varying output feedback law $w(t+1) = g(t, y(t), w(t))$, $u(t) = k(t, y(t), w(t))$, $w(t) \in W$ (Dynamic ROFS problem)

where

$$y(t) = h(t, x(t)), y(t) \in \mathcal{Y}'$$
(1.4)

is the so-called measured output, and $\mathcal{Y}', W$ are normed linear spaces.

In Section 4 of the present paper we present necessary and sufficient conditions for the solvability of the static and dynamic ROFS problem. To this end, we extend the notion of Complete Observability, introduced in [26] for continuous time systems (see also [3]). We remark that similar observability notions were also used in [1] for the construction of neural state estimators. Lyapunov-like conditions for the local solvability of the ROFS problem when the stabilized output is the whole state vector were given in [27] for autonomous finite-dimensional discrete-time systems without disturbances. The same problem was further studied in [2,8,14,19,20], where local and semi-global results were obtained. In the present paper it is shown that, if system (1.2) is stabilized by a continuous state feedback law and the feedback law is robustly completely observable from the measured output (1.4), then the dynamic ROFS problem for (1.2) is solvable. The procedure for the construction of the dynamic output feedback used in the proof of this result can be directly used for design purposes.

**Notation**
* By $\| \ \|_\mathcal{X}$, we denote the norm of the normed linear space $\mathcal{X}$. By $|\ |$ we denote the euclidean norm of $\Re^n$.
* $Z^+$ denotes the set of non-negative integers.
* For definitions of classes $K$, $K_\infty$, $KL$ see [17]. By $K^+$ we denote the set of all continuous positive functions defined on $\Re^+ := [0,+\infty)$.
* By $M_D$ we denote the set of all sequences $d = (d(0), d(1), d(2),....)$ with values in $D$, i.e., $d(i) \in D$ for all $i \in Z^+$.
* Let $H: Z^+ \times \mathcal{X} \to \mathcal{Y}$ a continuous map. The set-valued map $(t, y) \in Z^+ \times \mathcal{Y} \to H^{-1}(t, y) \subseteq \mathcal{X}$ is defined by $H^{-1}(t, y) = \{ x \in \mathcal{X}; H(t, x) = y \}$.
* Let $\mathcal{X}, \mathcal{Y}$ a pair of normed linear spaces. We denote by $CU(Z^+ \times A; W)$, where $A \subseteq \mathcal{X}$, the set of all continuous mappings $H: Z^+ \times A \to W \subseteq \mathcal{Y}$, with the following property: "for every pair of bounded sets $I \subset Z^+$, $S \subseteq A$ and for every $\varepsilon > 0$ the set $H(I \times S)$ is bounded and there exists $\delta > 0$ such that $\|H(t, x) - H(t, x_0)\|_\mathcal{Y} < \varepsilon$, for all $t \in I$, $x, x_0 \in S$ with $\|x - x_0\|_\mathcal{X} < \delta$".

## 2. Non-Uniform in Time Robust Global Asymptotic Output Stability (RGAOS)

In this section we first introduce the reader to the notion of non-uniform in time RGAOS for discrete-time systems as a special case of the notion of non-uniform in time RGAOS given in [11] for a wide class of systems. We consider the time-varying case (1.1) under the following hypotheses:

**(H1)** *There exist functions* $a \in K_\infty$, $\beta \in K^+$ *such that* $\|f(t,d,x)\|_X \le a(\beta(t)\|x\|_X)$, *for all* $(t,x,d) \in Z^+ \times X \times D$

**(H2):** *For every pair of bounded sets* $I \subset Z^+$, $S \subset X$ *and for every* $\varepsilon > 0$ *the set* $H(I \times S)$ *is bounded and there exists* $\delta > 0$ *such that* $\|H(t,x) - H(t,x_0)\|_Y < \varepsilon$, *for all* $t \in I$, $x, x_0 \in S$ *with* $\|x - x_0\|_X < \delta$. *Moreover, it holds that* $H(t,0) = 0$ *for all* $t \in Z^+$.

We note the following important fact for the time-varying case (1.1):

**Fact I:** *System (1.1) under hypothesis (H1) is Robustly Forward Complete (RFC) and there exist functions* $\mu \in K^+$, $a \in K_\infty$, *such that for every* $d \in M_D$, $(t_0, x_0) \in Z^+ \times X$, *the unique solution* $x(t)$ *of (1.1) initiated from* $x_0 \in X$ *at time* $t_0 \ge 0$ *and corresponding to* $d \in M_D$, *satisfies the following estimate:*

$$\|x(t)\|_X \le \mu(t)\, a(\|x_0\|_X), \quad \forall t \ge t_0 \tag{2.1}$$

Concerning the proof of Fact I, we notice that by virtue of Lemma 3.5 in [11] it suffices to show that system (1.1) under hypothesis (H1) is RFC and $0 \in X$ is a robust equilibrium point in the sense defined in [11]. Particularly, this follows by considering arbitrary $r \ge 0$, $T \in Z^+$, then defining recursively the sequence of sets in $X$ by $A(k) := f([0,2T] \times D \times A(k-1))$ for $k = 1,\ldots,T$ with $A(0) := \{x \in X; \|x\|_X \le r\}$, which are bounded by virtue of hypothesis (H1) and finally noticing that

$$\{x(t_0 + k, t_0, x_0; d)\,;\ \|x_0\|_X \le r,\ t_0 \le T,\ k \le T,\ d \in M_D\} \subseteq A(k) \text{ for all } k = 0,\ldots,T$$

where $x(t,t_0,x_0;d)$ denotes the unique solution of (1.1) initiated from $x_0 \in X$ at time $t_0 \ge 0$ and corresponding to $d \in M_D$. The fact that $0 \in X$ is a robust equilibrium point in the sense defined in [11] is an immediate consequence of hypothesis (H1) (details are left to the reader).

We are now in a position to present the definition of non-uniform in time RGAOS for discrete-time systems.

**Definition 2.1** *Let* $x(t)$ *denote the unique solution of (1.1) initiated from* $x_0 \in X$ *at time* $t_0 \ge 0$ *and corresponding to* $d \in M_D$. *We say that system (1.1) under hypotheses (H1-2) is **non-uniformly in time Robustly Globally Asymptotically Output Stable (RGAOS)** if it satisfies the following properties:*

**P1(Robust Output Stability)** *For every* $\varepsilon > 0$, $T \in Z^+$, *it holds that*

$$\sup\{\|H(t,x(t))\|_Y\,;\ t \ge t_0,\ \|x_0\|_X \le \varepsilon,\ t_0 \in [0,T],\ d \in M_D\} < +\infty$$

*and there exists a* $\delta := \delta(\varepsilon, T) > 0$ *such that:*

$$\|x_0\|_X \le \delta,\ t_0 \in [0,T] \Rightarrow \|H(t,x(t))\|_Y \le \varepsilon,\ \forall t \ge t_0,\ \forall d \in M_D$$

**P2(Robust Output Attractivity)** *For every* $\varepsilon > 0$, $T \in Z^+$ *and* $R \ge 0$, *there exists a* $\tau := \tau(\varepsilon, T, R) \in Z^+$, *such that:*

$$\|x_0\|_X \le R,\ t_0 \in [0,T] \Rightarrow \|H(t,x(t))\|_Y \le \varepsilon,\ \forall t \ge t_0 + \tau,\ \forall d \in M_D$$

*We say that system (1.1) is **non-uniformly in time strongly Robustly Globally Asymptotically Output Stable (strongly RGAOS)** if it is non-uniformly in time RGAOS and the set $H^{-1}(t,0) := \{ x \in X \, ; \, H(t,x) = 0 \}$ is positively invariant, i.e., if $x \in H^{-1}(t,0)$ then $H(t+1, f(t,d,x)) = 0$ for all $d \in D$.*

*Moreover, if there exists a function $a \in K_\infty$ such that $a(\|x\|_X) \leq \|H(t,x)\|_Y$ for all $(t,x) \in Z^+ \times X$, then we say that the equilibrium point $0 \in X$ is **non-uniformly in time Robustly Globally Asymptotically Stable (RGAS)** for system (1.1).*

The following facts are given in [11] as Lemma 3.3 and Lemma 3.4, respectively.

**Fact II:** *Suppose that (1.1) under hypotheses (H1-2) satisfies the Robust Output Attractivity Property (property P2 of Definition 2.1). Then (1.1) is non-uniformly in time RGAOS.*

**Fact III:** *System (1.1) under hypotheses (H1-2) is non-uniformly in time RGAOS if and only if there exist functions $\sigma \in KL$, $\beta \in K^+$ such that for every $d \in M_D$, $(t_0, x_0) \in Z^+ \times X$, the unique solution $x(t)$ of (1.1) initiated from $x_0 \in X$ at time $t_0 \geq 0$ and corresponding to $d \in M_D$, satisfies:*

$$\|H(t, x(t))\|_Y \leq \sigma\big(\beta(t_0)\|x_0\|_X, t-t_0\big), \ \forall t \geq t_0 \tag{2.2}$$

The following proposition provides Lyapunov-like characterizations for a time-varying system, which is non-uniformly in time RGAOS. It deals with the time-varying case (1.1) under the pair of hypotheses (H1-2) or under the pair of hypothesis (H2) and the following hypothesis:

**(H3):** *For every bounded sets $S \subset X$, $I \subset Z^+$ and for every $\varepsilon > 0$ the set $f(I \times D \times S)$ is bounded and there exists $\delta > 0$ such that $\sup\{\|f(t,d,x) - f(t,d,x_0)\|_X \, ; d \in D\} < \varepsilon$, for all $t \in I$, $x \in S$, $x_0 \in S$ with $\|x - x_0\|_X < \delta$. Moreover, it holds that $f(t,d,0) = 0$ for all $(t,d) \in Z^+ \times D$.*

**Remark about Hypothesis (H3):** Hypothesis (H3) is "stronger" hypothesis than (H1), in the sense that the implication (H3) $\Rightarrow$ (H1) holds. The proof of the implication (H3) $\Rightarrow$ (H1) is made by defining the following function:

$$a(T,s) := \sup\{\|f(t,d,x)\|_X \, ; t \in Z^+, t \leq T, \|x\|_X \leq s, d \in D\}$$

which is well-defined for all $T, s \geq 0$. Moreover, for every $T, s \geq 0$ the functions $a(\cdot, s)$ and $a(T, \cdot)$ are non-decreasing and since $f(t,0,d) = 0 \in X$ for all $(t,d) \in Z^+ \times D$, we also obtain $a(T,0) = 0$ for all $T \geq 0$. Finally, let $\varepsilon > 0$ and $T \geq 0$. It can be shown that hypothesis (H3) guarantees the existence of $\delta := \delta(\varepsilon, T) > 0$ such that $a(T, \delta(\varepsilon, T)) < \varepsilon$ and consequently we have $\lim_{s \to 0^+} a(T,s) = 0$ for all $T \geq 0$. It turns out from Lemma 2.3 in [10] that there exist functions $\zeta \in K_\infty$, $\beta \in K^+$ such that $a(T,s) \leq \zeta(\beta(T)s)$, for all $T, s \geq 0$ and consequently, hypothesis (H1) is satisfied.

**Proposition 2.2** *Consider system (1.1) under hypotheses (H1-2). Then the following statements are equivalent:*

**(i)** *System (1.1) is non-uniformly in time RGAOS.*

**(ii)** *There exist functions $V : Z^+ \times X \to \Re^+$, $a_1, a_2 \in K_\infty$, $\beta, \mu \in K^+$ and a constant $\lambda \in (0,1)$ such that:*

$$a_1\big(\|H(t,x)\|_Y + \mu(t)\|x\|_X\big) \leq V(t,x) \leq a_2\big(\beta(t)\|x\|_X\big), \ \forall (t,x) \in Z^+ \times X \tag{2.3a}$$

$$V(t+1, f(t,d,x)) \leq \lambda V(t,x), \ \forall (t,x,d) \in Z^+ \times X \times D \tag{2.3b}$$

*Moreover, if hypothesis (H3) holds then $V \in CU(Z^+ \times X; \Re^+)$.*

**(iii)** *There exist functions* $V: Z^+ \times \mathcal{X} \to \mathfrak{R}^+$, $a_1, a_2, a_3 \in K_\infty$ *with* $a_3(s) \leq s$ *for all* $s \geq 0$, $\beta \in K^+$ *and* $q \in C^0(\mathfrak{R}^+; \mathfrak{R}^+)$ *with* $\lim_{t \to +\infty} q(t) = 0$ *such that:*

$$a_1\left(\|H(t,x)\|_Y\right) \leq V(t,x) \leq a_2\left(\beta(t)\|x\|_\mathcal{X}\right), \quad \forall (t,x) \in Z^+ \times \mathcal{X} \tag{2.4a}$$

$$V(t+1, f(t,d,x)) \leq V(t,x) - a_3(V(t,x)) + q(t), \quad \forall (t,x,d) \in Z^+ \times \mathcal{X} \times D \tag{2.4b}$$

**Proof:** (i) $\Rightarrow$ (ii) We show the existence of a function $V: Z^+ \times \mathcal{X} \to \mathfrak{R}^+$, satisfying (2.3a,b), under the assumption of non-uniform in time RGAOS for (1.1). Since (1.1) under hypotheses (H1-2) is RGAOS, by virtue of Facts I and III, there exist functions $\mu \in K^+$, $a \in K_\infty$, $\sigma \in KL$, $\beta \in K^+$ such that for every $d \in M_D$, $(t_0, x_0) \in Z^+ \times \mathcal{X}$, the unique solution $x(t)$ of (1.1) initiated from $x_0 \in \mathcal{X}$ at time $t_0 \geq 0$ and corresponding to $d \in M_D$, satisfies (2.1) and (2.2). Next we consider the system:

$$\begin{aligned} z(t+1) &= \frac{\exp(-t-1)}{\mu(t+1)} f(t, d(t), \exp(t)\mu(t)z(t)) \\ w(t+1) &= \exp(-1)w(t) - \exp(-1)H(t, \exp(t)\mu(t)z(t)) + H(t+1, f(t,d(t),\exp(t)\mu(t)z(t))) \\ z(t) &\in \mathcal{X}, \; w(t) \in Y, \; t \in Z^+, \; d(t) \in D \end{aligned} \tag{2.5}$$

where the state space is the normed space $C := \mathcal{X} \times Y$ with norm $\|(z,w)\|_C := \sqrt{\|z\|_\mathcal{X}^2 + \|w\|_Y^2}$. We claim that zero for the above system is non-uniformly in time RGAS. Notice that the solution $(z(t), w(t))$ of (2.5) initiated from $(z_0, w_0) \in \mathcal{X} \times Y$ at time $t_0 \in Z^+$ and corresponding to $d(\cdot) \in M_D$ satisfies:

$$w(t) = H(t, \exp(t)\mu(t)z(t)) + \exp(-(t-t_0))(w_0 - H(t_0, \exp(t_0)\mu(t_0)z_0)), \quad \forall t \geq t_0 \tag{2.6}$$

Moreover the component $z(t)$ of the solution $(z(t), w(t))$ of (2.5) initiated from $(z_0, w_0) \in \mathcal{X} \times Y$ at time $t_0 \geq 0$ and corresponding to $d(\cdot) \in M_D$ is related to the solution $x(t)$ of (1.1) initiated from $x_0 = \exp(t_0)\mu(t_0)z_0$ at time $t_0 \in Z^+$ and corresponding to the same $d(\cdot) \in M_D$ with the following way:

$$x(t) = \exp(t)\mu(t)z(t), \quad \forall t \geq t_0 \tag{2.7}$$

Using (2.7) in conjunction with (2.1) and (2.2) we obtain:

$$\|z(t)\|_\mathcal{X} \leq \exp(-t) a\left(\exp(t_0)\mu(t_0)\|z_0\|_\mathcal{X}\right), \quad \forall t \geq t_0 \tag{2.8}$$

$$\|H(t, \exp(t)\mu(t)z(t))\|_Y \leq \sigma\left(\beta(t_0)\exp(t_0)\mu(t_0)\|z_0\|_\mathcal{X}, t-t_0\right), \quad \forall t \geq t_0 \tag{2.9}$$

Since $H(\cdot)$ is continuous with $H(t,0) = 0$ for all $t \geq 0$ and since the set $H(I \times S)$ is bounded for every pair of bounded sets $I \subset Z^+$ and $S \subset \mathcal{X}$, it follows from Lemma 3.2 in [11] that there exist functions $\zeta \in K_\infty$ and $\gamma \in K^+$ such that

$$\|H(t,x)\|_Y \leq \zeta(\gamma(t)\|x\|_\mathcal{X}), \quad \forall (t,x) \in Z^+ \times \mathcal{X} \tag{2.10}$$

Combining estimate (2.9) with (2.6) and inequality (2.10), we obtain for all $t \geq t_0$:

$$\|w(t)\|_Y \leq \sigma\left(\beta(t_0)\exp(t_0)\mu(t_0)\|z_0\|_\mathcal{X}, t-t_0\right) + \exp(-(t-t_0))\left[\|w_0\|_Y + \zeta(\gamma(t_0)\exp(t_0)\mu(t_0)\|z_0\|_\mathcal{X})\right] \tag{2.11}$$

We conclude from (2.8) and (2.11) that the solution $(z(t), w(t))$ of (2.5) initiated from $(z_0, w_0) \in \mathcal{X} \times Y$ at time $t_0 \in Z^+$ and corresponding to $d(\cdot) \in M_D$ satisfies:

$$\|(z(t), w(t))\|_C \leq \tilde{\sigma}\left(\tilde{\beta}(t_0) \|(z_0, w_0)\|_C, t - t_0\right), \quad \forall t \geq t_0 \tag{2.12}$$

where $\tilde{\sigma}(s,t) := a(s)\exp(-t) + \sigma(s,t) + \exp(-t)(s + \zeta(s))$, $\tilde{\beta}(t) := 1 + (1 + \beta(t) + \gamma(t))\exp(t)\mu(t)$. Inequality (2.12) implies that zero is non-uniformly in time RGAS for (2.5). Moreover, since the dynamics of system (2.5)

$$\tilde{f}(t,d,z,w) := \begin{pmatrix} \dfrac{\exp(-t-1)}{\mu(t+1)} f(t,d,\exp(t)\mu(t)z) \\ \exp(-1)w - \exp(-1)H(t,\exp(t)\mu(t)z) + H(t+1, f(t,d,\exp(t)\mu(t)z)) \end{pmatrix}$$ 

satisfy hypothesis (H1), it follows from Theorem 2.9 in [12] that there exist functions $U : Z^+ \times X \times Y \to \Re^+$, $\tilde{a}_1(\cdot), \tilde{a}_2(\cdot) \in K_\infty$, $\beta_2(\cdot) \in K^+$ such that:

$$\tilde{a}_1\left(\|(z,w)\|_C\right) \leq U(t,z,w) \leq \tilde{a}_2\left(\beta_2(t)\|(z,w)\|_C\right), \quad \forall (t,x,w) \in Z^+ \times C \tag{2.13a}$$

$$U\left(t+1, \dfrac{\exp(-t-1)}{\mu(t+1)} f(t,d,\exp(t)\mu(t)z), \exp(-1)(w - H(t,\exp(t)\mu(t)z)) + H(t+1, f(t,d,\exp(t)\mu(t)z))\right) \leq \lambda U(t,z,w)$$

$$\forall (t,z,w,d) \in Z^+ \times C \times D \tag{2.13b}$$

Finally, we define:

$$V(t,x) := U\left(t, \dfrac{\exp(-t)}{\mu(t)} x, H(t,x)\right) \tag{2.14}$$

Inequality (2.3b) is an immediate consequence of (2.13b) and definition (2.14). Moreover, by virtue of Lemma 2.3 in [10], there exist functions $a_2 \in K_\infty$ and $\beta \in K^+$ such that

$$\tilde{a}_2\left(\dfrac{\beta_2(t)\exp(-t)}{\mu(t)} s + \beta_2(t)\zeta(\gamma(t)s)\right) \leq a_2(\beta(t)s), \quad \forall (t,s) \in \Re^+ \times \Re^+ \tag{2.15}$$

Inequalities (2.13a), (2.10), (2.15) and the trivial inequality $\max(\|w\|_Y, \|z\|_X) \leq \|(z,w)\|_C \leq \|z\|_X + \|w\|_Y$ imply inequality (2.3a) with $\mu(t) := \dfrac{\exp(-t)}{\mu(t)}$ and $a_1(s) := \tilde{a}_1\left(\dfrac{s}{2}\right)$. Finally, if hypothesis (H3) is satisfied for (1.1), it follows that the dynamics of system (2.5)

$$\tilde{f}(t,d,z,w) := \begin{pmatrix} \dfrac{\exp(-t-1)}{\mu(t+1)} f(t,d,\exp(t)\mu(t)z) \\ \exp(-1)w - \exp(-1)H(t,\exp(t)\mu(t)z) + H(t+1, f(t,d,\exp(t)\mu(t)z)) \end{pmatrix}$$ also satisfy hypothesis (H3).

Consequently, by virtue of Theorem 2.9 in [12], we conclude that $U \in CU(Z^+ \times X \times Y; \Re^+)$ and since $H \in CU(Z^+ \times X; Y)$ we obtain that $V \in CU(Z^+ \times X; \Re^+)$.

(ii) $\Rightarrow$ (iii) This implication is trivial (notice that statement (iii) holds with $q(t) \equiv 0$ and $a_3(s) := (1-\lambda)s$).

(iii) $\Rightarrow$ (i) Assuming that statement (iii) holds, we will show that system (1.1) satisfies the Robust Output Attractivity Property (property P2 of Definition 2.1). Then, by virtue of Fact II, statement (i) holds. Without loss of generality we may assume that $q(t) > 0$ for all $t \in Z^+$. Let arbitrary $\varepsilon > 0$, $T \in Z^+$, $R \geq 0$, $\|x_0\|_X \leq R$, $t_0 \in [0,T]$ and $d \in M_D$. Let also $V(t) := V(t,x(t))$ and $q_{\max} := \sup_{t \geq 0} q(t)$, where $x(t)$ is the unique solution of (1.1) corresponding to $d \in M_D$, initiated from $x_0$ at time $t_0$. It follows from inequality (2.4b) and Lemma 3.1 in [12] that:

$$V(t) \leq V(t_0) + a_3^{-1}(q_{\max}) + q_{\max}, \quad \forall t \geq t_0 \tag{2.16}$$

$$V(t) < a_3^{-1}\left(2 \sup_{t \geq t_0 + \tilde{\tau}} q(t)\right) + \sup_{t \geq t_0 + \tilde{\tau}} q(t), \quad \forall \tilde{\tau} \in Z^+, \forall t \geq t_0 + \tilde{\tau} + \dfrac{V(t_0 + \tilde{\tau})}{\sup_{t \geq t_0 + \tilde{\tau}} q(t)} \tag{2.17}$$

Since $\lim_{t \to +\infty} q(t) = 0$, there exists $\tilde{\tau} := \tilde{\tau}(\varepsilon) \in Z^+$ such that $a_3^{-1}(2 \sup_{t \geq \tilde{\tau}} q(t)) + \sup_{t \geq \tilde{\tau}} q(t) \leq a_1(\varepsilon)$. Combining the latter inequality with (2.4a), (2.16), (2.17), we conclude that the Robust Output Attractivity Property holds for system (1.1) with

$$\tau(\varepsilon, T, R) := T + \tilde{\tau}(\varepsilon) + \left[ \frac{a_2\left(R \max_{0 \leq t_0 \leq T} \beta(t_0)\right) + a_3^{-1}(q_{max}) + q_{max}}{\sup_{t \geq T+\tilde{\tau}(\varepsilon)} q(t)} \right] + 1$$

where $[x]$ denotes the integer part of the real number $x$. The proof is complete. ◁

**Example 2.3** Consider the nonlinear finite-dimensional discrete-time time-varying system:

$$\begin{aligned}
x_1(t+1) &= d(t) x_1(t) \\
x_2(t+1) &= 2^{-t} d(t) |x_1(t)|^{\frac{1}{2}} \\
Y(t) &= H(t, x(t)) := x_2(t) \\
x(t) &:= (x_1(t), x_2(t)) \in \Re^2, \, t \in Z^+, \, d(t) \in [-2,2]
\end{aligned} \quad (2.18)$$

Consider the continuous function $V(t,x) := \exp(-t)|x_1| + |x_2|$, which clearly satisfies the following inequality:

$$|Y| = |x_2| \leq V(t,x) \leq 2|x|, \, \forall (t,x) \in Z^+ \times \Re^2 \quad (2.19)$$

Moreover, notice that for all $(t, x, d) \in Z^+ \times \Re^2 \times [-2,2]$ we obtain:

$$V(t+1, dx_1, 2^{-t} d|x_1|^{\frac{1}{2}}) = \exp(-t-1)|d||x_1| + 2^{-t}|d||x_1|^{\frac{1}{2}} \leq 2e^{-1} \exp(-t)|x_1| + 2^{-t+1}|x_1|^{\frac{1}{2}}$$

$$\leq \frac{2+e}{2e} \exp(-t)|x_1| + \frac{2e}{e-2}\left(\frac{e}{4}\right)^t \leq \lambda V(t,x) + q(t) \quad (2.20)$$

where $\lambda := \frac{2+e}{2e} \in (0,1)$ and $q(t) := \frac{2e}{e-2}\left(\frac{e}{4}\right)^t$ with $\lim_{t \to +\infty} q(t) = 0$. By virtue of (2.19) and (2.20) it follows that statement (iii) of Proposition 2.2 is satisfied with $\beta(t) \equiv 2$, $a_1(s) = a_2(s) := s$ and $a_3(s) := (1-\lambda)s$. We conclude that system (2.18) is non-uniformly in time RGAOS. ◁

## 3. Non-Uniform in Time Input-to-Output Stability (IOS)

In this section we first introduce the reader to the notion of non-uniform in time IOS for discrete-time systems as a special case of the notion of non-uniform in time IOS given in [11] for a wide class of systems. We consider the time-varying case (1.2) under hypothesis (H2) and the following hypothesis:

**(A1)** *There exist functions* $a \in K_\infty$, $\beta \in K^+$ *such that* $\|f(t,d,x,u)\|_X \leq a\left(\beta(t)\|x\|_X\right) + a\left(\beta(t)\|u\|_U\right)$, *for all* $(t, x, d, u) \in Z^+ \times X \times D \times U$.

First we note the following important fact for the time-varying case (1.2):

**Fact IV:** *System (1.2) under hypothesis (A1) is Robustly Forward Complete (RFC) from the input* $u \in M_U$ *and there exist functions* $\mu \in K^+$, $a \in K_\infty$ *and a constant* $R \geq 0$ *such that for every* $(t_0, x_0, d, u) \in Z^+ \times X \times M_D \times M_U$, *the corresponding solution* $x(t)$ *of (1.2) with* $x(t_0) = x_0$ *satisfies the following estimate:*

$$\|x(t)\|_X \leq \mu(t) a\left(R + \|x_0\|_X + \sup_{\tau \in [t_0, t]} \|u(\tau)\|_U\right), \, \forall t \geq t_0 \quad (3.1)$$

Concerning the proof of Fact IV, we notice that by virtue of Lemma 3.5 in [11] it suffices to show that system (1.2) under hypothesis (A1) is RFC from the input $u \in M_U$. Particularly, this follows by considering arbitrary $r \geq 0$, $T \in Z^+$, then defining recursively the sequence of sets in $X$ by $A(k) := f([0, 2T] \times D \times A(k-1) \times B[0, r])$ for $k = 1, ..., T$, where $B[0, r] := \{u \in U ; \|u\|_U \leq r\}$ and $A(0) := \{x \in X ; \|x\|_X \leq r\}$, which are bounded by virtue of hypothesis (A1) and finally noticing that

$$\{x(t_0+k, t_0, x_0; u, d) ; \|x_0\|_X \leq r, t_0 \leq T, k \leq T, d \in M_D, u \in M_{B[0,r]}\} \subseteq A(k) \text{ for all } k = 0, ..., T$$

where $x(t, t_0, x_0; u, d)$ denotes the unique solution of (1.1) initiated from $x_0 \in X$ at time $t_0 \geq 0$ and corresponding to input $(u, d) \in M_{B[0,r]} \times M_D$.

We are now in a position to present the definition of non-uniform in time IOS property for discrete-time systems.

**Definition 3.1:** *Let $x(t)$ denote the unique solution of (1.2) initiated from $x_0 \in X$ at time $t_0 \in Z^+$ and corresponding to $(d, u) \in M_D \times M_U$. We say that system (1.2) under the pair of hypotheses (A1), (H2) satisfies the **non-uniform in time Input-to-Output Stability property (IOS)** from the input $u \in M_U$ if there exist functions $\sigma \in KL$, $\beta, \gamma \in K^+$ and $\rho \in K_\infty$ such that the following estimate holds for all $u \in M_U$, $(t_0, x_0, d) \in Z^+ \times X \times M_D$ and $t \in [t_0, +\infty)$:*

$$\|H(t, x(t))\|_Y \leq \max\left\{\sigma(\beta(t_0)\|x_0\|_X, t-t_0), \sup_{\tau \in [t_0, t]} \sigma(\beta(\tau)\rho(\gamma(\tau)\|u(\tau)\|_U), t-\tau)\right\} \quad (3.2)$$

*Moreover, if there exists a function $a \in K_\infty$ such that $a(\|x\|_X) \leq \|H(t, x)\|_Y$ for all $(t, x) \in Z^+ \times X$, then we say that (1.2) satisfies the **non-uniform in time Input-to-State Stability property (ISS)** from the input $u \in M_U$.*

The following proposition provides various characterizations of the non-uniform in time IOS property for the time-varying case (1.2). Moreover, it shows how the notion of non-uniform in time IOS is related to the notion of non-uniform in time RGAOS. It deals with the time-varying case (1.2) under the pair of hypotheses (A1) and (H2) or under the pair of hypothesis (H2) and the following hypothesis:

**(A2):** *For every bounded sets $S \subset X \times U$, $I \subset Z^+$ and for every $\varepsilon > 0$, the set $f(I \times D \times S)$ is bounded and there exists $\delta > 0$ such that $\sup\{\|f(t, d, x, u) - f(t, d, x_0, u_0)\|_X ; d \in D\} < \varepsilon$, for all $t \in I$, $(x, u) \in S$, $(x_0, u_0) \in S$ with $\|x - x_0\|_X + \|u - u_0\|_U < \delta$. Moreover, it holds that $f(t, d, 0, 0) = 0$ for all $(t, d) \in Z^+ \times D$.*

**Remark about Hypothesis (A2):** Hypothesis (A2) is "stronger" hypothesis than (A1), in the sense that the implication (A2) $\Rightarrow$ (A1) holds. The proof of the implication (A2) $\Rightarrow$ (A1) is made by defining the following function:

$$a(T, s) := \sup\{\|f(t, d, x, u)\|_X ; t \in Z^+, t \leq T, \|x\|_X \leq s, d \in D, \|u\|_U \leq s\}$$

which is well-defined for all $T, s \geq 0$. Moreover, for every $T, s \geq 0$ the functions $a(\cdot, s)$ and $a(T, \cdot)$ are non-decreasing and since $f(t, 0, d, 0) = 0 \in X$ for all $(t, d) \in Z^+ \times D$, we also obtain $a(T, 0) = 0$ for all $T \geq 0$. Finally, let $\varepsilon > 0$ and $T \geq 0$. It can be shown that hypothesis (A2) guarantees the existence of $\delta := \delta(\varepsilon, T) > 0$ such that $a(T, \delta(\varepsilon, T)) < \varepsilon$ and consequently we have $\lim_{s \to 0^+} a(T, s) = 0$ for all $T \geq 0$. It turns out from Lemma 2.3 in [10] that there exist functions $\zeta \in K_\infty$, $\beta \in K^+$ such that $a(T, s) \leq \zeta(\beta(T)s)$, for all $T, s \geq 0$. Consequently, we obtain $\|f(t, x, d, u)\|_X \leq \zeta(\beta(t) \max\{\|x\|_X, \|u\|_U\})$, for all $(t, x, d, u) \in Z^+ \times X \times D \times U$, which directly implies hypothesis (A1).

**Proposition 3.2** *Consider system (1.2) under the pair of hypotheses (A1), (H2). Then the following statements are equivalent:*

**(i)** *System (1.2) satisfies the non-uniform in time IOS property.*

**(ii)** *There exist functions $\zeta \in K_\infty$, $\beta, \delta \in K^+$ and $\sigma \in KL$ such that for every $(t_0, x_0, d, u) \in Z^+ \times X \times M_D \times M_U$, the corresponding solution $x(t)$ of (1.2) with $x(t_0) = x_0$ satisfies:*

$$\|H(t, x(t))\|_Y \leq \max\left\{\sigma(\beta(t_0)\|x_0\|_X, t - t_0), \sup_{t_0 \leq \tau \leq t} \zeta(\delta(\tau)\|u(\tau)\|_U)\right\}, \quad \forall t \geq t_0 \tag{3.3}$$

**(iii)** *There exist functions $\theta \in K_\infty$ and $p \in K^+$ such that the following system is non-uniformly in time RGAOS:*

$$\begin{aligned} x(t+1) &= f(t, d(t), x(t), p(t)\theta(\|x(t)\|_X)d'(t)) \\ Y(t) &= H(t, x(t)) \\ x(t) &\in X, \ Y(t) \in Y, \ (d(t), d'(t)) \in D \times B[0,1], \ t \in Z^+ \end{aligned} \tag{3.4}$$

*where $B[0,1] := \{u \in U; \|u\|_U \leq 1\}$.*

**(iv)** *There exist functions $V : Z^+ \times X \to \Re^+$, $a_1, a_2, a_3 \in K_\infty$, $\beta, \phi, \mu \in K^+$ and a constant $\lambda \in (0,1)$ such that:*

$$a_1(\|H(t,x)\|_Y + \mu(t)\|x\|_X) \leq V(t,x) \leq a_2(\beta(t)\|x\|_X), \quad \forall (t,x) \in Z^+ \times X \tag{3.5a}$$

$$V(t+1, f(t,d,x,u)) \leq \lambda V(t,x) + a_3(\phi(t)\|u\|_U), \quad \forall (t,x,d,u) \in Z^+ \times X \times D \times U \tag{3.5b}$$

*Moreover, if hypothesis (A2) holds then $V \in CU(Z^+ \times X; \Re^+)$.*

**Proof:** (i)$\Rightarrow$(ii) Suppose that (1.2) satisfies the non-uniform in time IOS property. Then there exist functions $\sigma \in KL$, $\beta, \mu, \gamma \in K^+$, $a \in K_\infty$ and $\rho \in K_\infty$ such that (3.2) holds for all $u \in M_U$, $(t_0, x_0, d) \in Z^+ \times X \times M_D$ and $t \in [t_0, +\infty)$. By invoking Lemma 2.3 in [10], there exist functions $a \in K_\infty$ and $\delta \in K^+$ such that $\beta(t)\rho(\phi(t)s) \leq a(\delta(t)s)$ for all $(t,s) \in \Re^+ \times \Re^+$ and if we set $\zeta(s) := \sigma(a(s), 0) + s$ (which obviously is of class $K_\infty$), the desired (3.3) is a consequence of (3.2) and the previous inequality.

(ii)$\Rightarrow$(iii) Suppose that there exist functions $\zeta \in K_\infty$, $\beta, \delta \in K^+$ and $\sigma \in KL$ such that for every $(t_0, x_0, d, u) \in Z^+ \times X \times M_D \times M_U$, the corresponding solution $x(t)$ of (1.2) with $x(t_0) = x_0$ satisfies (3.3). Moreover, by virtue of Fact IV, there exist functions $\mu \in K^+$, $a \in K_\infty$ and a constant $R \geq 0$ such that for every $(t_0, x_0, d, u) \in Z^+ \times X \times M_D \times M_U$, the corresponding solution $x(t)$ of (1.2) with $x(t_0) = x_0$ satisfies (3.1). Consequently, for every $(t_0, x_0, d, u) \in Z^+ \times X \times M_D \times M_U$, the corresponding solution $x(t)$ of (1.2) with $x(t_0) = x_0$ satisfies the following estimate:

$$\frac{\|x(t)\|_X}{\mu(t)} \leq a(2R + 2\|x_0\|_X) + \sup_{\tau \in [t_0, t]} a(2\|u(\tau)\|_U), \quad \forall t \geq t_0 \tag{3.6}$$

Without loss of generality we may assume that the functions $\zeta \in K_\infty$ and $\delta \in K^+$ involved in (3.3) satisfy $\zeta(s) \geq s$ (or equivalently $\zeta^{-1}(s) \leq s$) and $\delta(t) \geq 1$ for all $t, s \geq 0$ and that the function $\mu \in K^+$ involved in (3.1) is non-decreasing. Define:

$$\gamma(t, s) := 2\mu(t)a(2(1 + \mu(t))\exp(t)\zeta(\delta(t)s)) \tag{3.7}$$

By virtue of Lemma 2.3 in [10] there exist functions $q \in K^+$, $a' \in K_\infty$ such that $\gamma(t,s) \leq a'(q(t)s)$ for all $t,s \geq 0$. Define:

$$\theta^{-1}(s) := a'(s) \text{ and } p(t) := \frac{1}{q(t)} \tag{3.8}$$

It follows from definitions (3.7) and (3.8) that:

$$p(t)\theta(s) \leq \frac{1}{\delta(t)} \zeta^{-1}\left(\frac{\exp(-t)}{2(1+\mu(t))} a^{-1}\left(\frac{s}{2\mu(t)}\right)\right), \text{ for all } t,s \geq 0 \tag{3.9}$$

Notice that since (1.2) satisfies the pair of hypotheses (A1) and (H2), it follows that system (3.4) satisfies hypotheses (H1-2). Clearly, the solution of system (3.4) with $x(t_0) = x_0$ corresponding to $(d,d') \in M_D \times M_{B[0,1]}$ coincides with the solution of (1.2) with same initial condition corresponding to inputs $d \in M_D$ and $u \in M_U$ with $u(t) = p(t)\theta(\|x(t)\|_X)d'(t)$ for all $t \geq t_0$. Consequently, since $\zeta^{-1}(s) \leq s$ and $\delta(t) \geq 1$ for all $t,s \geq 0$, we obtain from (3.9) that $\|u(t)\|_U \leq \frac{1}{2} a^{-1}\left(\frac{\|x(t)\|_X}{2\mu(t)}\right)$ for all $t \geq t_0$ and it follows from (3.6) that:

$$\sup_{\tau \in [t_0,t]}\left(\frac{\|x(\tau)\|_X}{\mu(\tau)}\right) \leq 2a(2R+2\|x_0\|_X), \forall t \geq t_0 \tag{3.10}$$

Combining inequalities (3.3) and (3.9) we obtain:

$$\|H(t,x(t))\|_Y \leq \max\left\{\sigma(\beta(t_0)\|x_0\|_X, t-t_0), \frac{\exp(-t_0)}{2(1+\mu(t_0))} \sup_{t_0 \leq \tau \leq t} a^{-1}\left(\frac{\|x(\tau)\|_X}{2\mu(\tau)}\right)\right\}, \forall t \geq t_0 \tag{3.11}$$

Estimate (3.11) in conjunction with estimate (3.10) gives:

$$\|H(t,x(t))\|_Y \leq \max\left\{\sigma(\beta(t_0)\|x_0\|_X, t-t_0), \exp(-t_0)\left(R + \frac{\|x_0\|_X}{1+\mu(t_0)}\right)\right\}, \forall t \geq t_0 \tag{3.12}$$

Notice that by virtue of (3.10) and (3.12) we obtain for all $t_1 \in Z^+$ and $t \geq t_0 + t_1$:

$$\|H(t,x(t))\|_Y \leq \max\left\{\sigma(\beta(t_0+t_1)\mu(t_0+t_1)2a(2R+2\|x_0\|_X), t-t_0-t_1), \exp(-t_0-t_1)(R+2a(2R+2\|x_0\|_X))\right\} \tag{3.13}$$

Next we establish robust global asymptotic output stability. Without loss of generality we may assume that the function $\beta \in K^+$ involved in (3.13) is non-decreasing. Consider the function $b(t,T,r) := \sup\{\|H(t_0+t, x(t_0+t))\|_Y ; (d,d') \in M_D \times M_{B[0,1]}, \|x_0\|_X \leq r, t_0 \in [0,T]\}$, where $x(\cdot)$ denotes the solution of (3.4) with $x(t_0) = x_0$ corresponding to some $(d,d') \in M_D \times M_{B[0,1]}$. It suffices to show that $\lim_{t \to +\infty} b(t,T,r) = 0$, for all $(T,r) \in (\Re^+)^2$. Let $\varepsilon > 0$ arbitrary. Clearly, there exists $t_1 = t_1(\varepsilon,r) \in Z^+$ and $t_2 = t_2(\varepsilon,T,r) \in Z^+$ such that $(R+2a(2R+2r))\exp(-t_1) \leq \varepsilon$ and $\sigma(\beta(T+t_1)\mu(T+t_1)2a(2R+2r), t_2) \leq \varepsilon$. By virtue of (3.13) and definition of $b$, we obtain for all $t \geq t_1$:

$$b(t,T,r) \leq \max\{\sigma(\beta(T+t_1)\mu(T+t_1)2a(2R+2r), t-t_1), \varepsilon\}$$

which directly implies that $b(t,T,r) \leq \varepsilon$, for all $t \geq t_1 + t_2$. Thus the robust output attractivity property is satisfied for $\tau(\varepsilon,T,r) := t_1(\varepsilon,r) + t_2(\varepsilon,T,r)$. By virtue of Fact II, we conclude that system (3.4) is non-uniformly in time RGAOS.

(iii) ⇒ (iv) Notice that since (1.2) satisfies the pair of hypotheses (A1) and (H2) (or (A2) and (H2)), it follows that system (3.4) satisfies hypotheses (H1-2) (or (H2-3)). Since (3.4) is non-uniformly in time RGAOS and satisfies hypotheses (H1-2) (or (H2-3)), it follows from Proposition 2.2, that there exist functions $V : Z^+ \times \mathcal{X} \to \Re^+$ ($V \in CU(Z^+ \times \mathcal{X}; \Re^+)$), $a_1$, $a_2$ of class $K_\infty$, $\beta, \mu$ of class $K^+$ and constant $\lambda \in (0,1)$ such that:

$$a_1\left(\|H(t,x)\|_Y + \mu(t)\|x\|_\mathcal{X}\right) \leq V(t,x) \leq a_2\left(\beta(t)\|x\|_\mathcal{X}\right), \quad \forall (t,x) \in Z^+ \times \mathcal{X} \quad (3.14a)$$

$$V(t+1, f(t,d,x,u)) \leq \lambda V(t,x), \quad \forall (t,x,d,u) \in Z^+ \times \mathcal{X} \times D \times U, \text{ with } \|u\|_U \leq p(t)\theta\left(\|x\|_\mathcal{X}\right) \quad (3.14b)$$

Define for all $(t,x,u) \in Z^+ \times \mathcal{X} \times U$:

$$\psi(t,x,u) := \sup\{V(t+1, f(t,d,x,u)) ; \ d \in D \} \quad (3.15)$$

Clearly, hypothesis (A1) (which holds in any case; see Remark about hypothesis (A2) above), inequality (3.14a) in conjunction with Lemma 2.3 in [10] imply the existence of functions $\omega \in K_\infty$ and $q \in K^+$ such that $\psi(t,x,u) \leq \omega(q(t)\|x\|_\mathcal{X}) + \omega(q(t)\|u\|_U)$. Moreover, Lemma 2.3 in [10] guarantees the existence of functions $a_3 \in K_\infty$ and $\phi \in K^+$ such that $\omega\left(q(t)\theta^{-1}\left(\dfrac{s}{p(t)}\right)\right) + \omega(q(t)s) \leq a_3(\phi(t)s)$ for all $t,s \geq 0$. Combining the previous inequalities and definition (3.15) we obtain:

$$\sup\left\{V(t+1, f(t,d,x,u)) ; \ d \in D, \|u\|_U \leq s, \|x\|_\mathcal{X} \leq \theta^{-1}\left(\dfrac{s}{p(t)}\right)\right\} \leq a_3(\phi(t)s), \text{ for all } t,s \geq 0 \quad (3.16)$$

We next establish inequality (3.5b), with $a_3$ as previously, by considering the following two cases:
* $\|u\|_U \leq p(t)\theta(\|x\|_\mathcal{X})$. In this case inequality (3.5b) is a direct consequence of (3.14b).
* $\|u\|_U \geq p(t)\theta(\|x\|_\mathcal{X})$. In this case inequality (3.5b) is a direct consequence of (3.16).

(iv) ⇒ (i) Consider the trajectory $x(t)$ of (1.2) that corresponds to input $(d,u) \in M_D \times M_U$ with initial condition $x(t_0) = x_0 \in \mathcal{X}$ and let $c := -\log(\lambda) > 0$, $V(t) = V(t, x(t))$, $b(t) := \exp(2ct)a_3(\phi(t)\|u(t)\|_U)$ for all $t \geq t_0$. Inequality (3.5b) implies that $V(t+1) \leq \exp(-c)V(t) + \exp(-2ct)b(t)$ for all $t \geq t_0$, which gives (using induction arguments):

$$V(t) \leq \exp(-c(t-t_0))V(t_0) + \dfrac{\exp(2c)}{\exp(c)-1}\exp(-c(t-t_0))\sup_{t_0 \leq \tau \leq t}\left(\exp(2c\tau)a_3(\phi(\tau)\|u(\tau)\|_U)\right), \text{ for all } t \geq t_0 \quad (3.17)$$

By Lemma 2.3 in [10] there exist functions $\rho \in K_\infty$ and $\gamma \in K^+$ such that $\dfrac{1}{\beta(t)}a_2^{-1}\left(\dfrac{\exp(2c)}{\exp(c)-1}\exp(2ct)a_3(\phi(t)s)\right) \leq \rho(\gamma(t)s)$, where $\beta \in K^+$, $a_2 \in K_\infty$ are the functions involved in (3.5a). The previous inequality, in conjunction with inequality (3.5a), (3.17) and definition $\sigma(s,t) := 2a_1^{-1}(2\exp(-ct)a_2(s)) \in KL$ implies (3.2). The proof is complete. ◁

The following proposition provides a sharper characterization of the IOS property for the time-varying case (1.2), which holds only for discrete-time systems with continuous dynamics. For continuous-time systems the situation is more involved since the finite escape time phenomenon can occur (see [24,25]). Further research is required for the case of discrete-time systems with discontinuous dynamics.

**Proposition 3.3** *System (1.2) under the pair of hypotheses (A2), (H2), satisfies the non-uniform in time IOS property if and only if the "unforced" system (1.3) is non-uniformly in time RGAOS.*

**Proof** It is clear that if system (1.2) satisfies the non-uniform in time IOS property then the "unforced" system (1.3) is non-uniformly in time RGAOS. Therefore it suffices to prove the converse statement. Specifically, we show that if the "unforced" system (1.3) is non-uniformly in time RGAOS, then statement (ii) of Proposition 3.2 holds for system (1.2). Consequently, by equivalence of statements (i) and (ii) of Proposition 3.2 it follows that system (1.2) satisfies the non-uniform in time IOS property.

Notice that by virtue of hypothesis (A2), the "unforced" system (1.3) satisfies hypothesis (H3). Since the "unforced" system (1.3) is non-uniformly in time RGAOS, it follows by Proposition 2.2, that there exist functions $V \in CU(Z^+ \times X ; \Re^+)$, $a_1$, $a_2$ of class $K_\infty$, $\tilde{\beta}$ of class $K^+$ and constant $c > 0$ such that:

$$a_1\left(\|H(t,x)\|_Y\right) \leq V(t,x) \leq a_2\left(\beta(t)\|x\|_X\right), \quad \forall (t,x) \in Z^+ \times X \tag{3.18a}$$

$$\sup_{d \in D} V(t+1, f(t,d,x,0)) \leq \exp(-c) V(t,x), \quad \forall (t,x,d) \in Z^+ \times X \times D \tag{3.18b}$$

Define the following function:

$$\gamma(r,s) := \sup\left\{ \left|V(t+1, f(t,d,x,u)) - V(t+1, f(t,d,x,0))\right|; 0 \leq t \leq r, d \in D, \|x\|_X \leq r, \|u\|_U \leq s \right\} \tag{3.19}$$

Clearly, by virtue of the right-hand side inequality (3.18a) and hypothesis (A1) (which holds since hypothesis (A2) holds; see Remark about Hypothesis (A2) above), it follows that $\gamma(r,s) < +\infty$ for all $r, s \geq 0$. Moreover, definition (3.19) guarantees that for every $r, s \geq 0$ the mappings $\gamma(r, \cdot)$ and $\gamma(\cdot, s)$ are non-decreasing with $\gamma(r,0) = 0$. Finally, hypothesis (A2) in conjunction with the fact that $V \in CU(Z^+ \times X ; \Re^+)$ guarantees that $\lim_{s \to 0^+} \gamma(r,s) = 0$ for all $r \geq 0$. Consequently, Lemma 2.3 in [10] guarantees the existence of functions $a_3 \in K_\infty$ and $\phi \in K^+$ such that $\gamma(r,s) \leq a_3(\phi(r)s)$ for all $r, s \geq 0$. It follows by definition (3.19) that we have for all $(t,x,u) \in Z^+ \times X \times U$:

$$\left|\sup_{d \in D} V(t+1, f(t,d,x,u)) - \sup_{d \in D} V(t+1, f(t,d,x,0))\right| \leq a_3\left(\phi(t)\|u\|_U\right) + a_3\left(\phi(\|x\|_X)\|u\|_U\right) \tag{3.20}$$

By virtue of Fact IV, there exist functions $\mu \in K^+$, $a \in K_\infty$ and a constant $R \geq 0$ such that for every $(t_0, x_0, d, u) \in Z^+ \times X \times M_D \times M_U$, the corresponding solution $x(t)$ of (1.2) with $x(t_0) = x_0$ satisfies:

$$\frac{\|x(t)\|_X}{\mu(t)} \leq a(4R) + a(4\|x_0\|_X) + \sup_{\tau \in [t_0, t]} a(2\|u(\tau)\|_U), \quad \forall t \geq t_0 \tag{3.21}$$

Using (3.18b), (3.20), Lemma 2.3 in [10] and Corollary 10 and Remark 11 in [23], we obtain functions $a_5, a_6 \in K_\infty$, $q \in K^+$ such that for all $(t,x,d,u) \in Z^+ \times X \times D \times U$ it holds that:

$$\sup_{d \in D} V(t+1, f(t,d,x,u)) \leq \exp(-c) V(t,x) + \exp(-2ct) a_5\left(q(t)\|u\|_U\right) + \exp(-2ct) a_6\left(\frac{\|x\|_X}{\mu(t)}\right) a_5\left(q(t)\|u\|_U\right) \tag{3.22}$$

Let $(t_0, x_0, d, u) \in Z^+ \times X \times M_D \times M_U$ and consider the corresponding solution $x(t)$ of (1.2) with $x(t_0) = x_0$. Let $V(t) := V(t, x(t))$ for $t \geq t_0$. By virtue of (3.21) and (3.22) we obtain:

$$V(t+1) \leq \exp(-c) V(t) + \exp(-2ct) \sup_{\tau \in [t_0, t]} a_5\left(q(\tau)\|u(\tau)\|_U\right) + \exp(-2ct) a_6(4a(4R)) \sup_{\tau \in [t_0, t]} a_5\left(q(\tau)\|u(\tau)\|_U\right)$$

$$+ \exp(-2ct) \sup_{\tau \in [t_0, t]} \left(a_5(q(\tau)\|u(\tau)\|_U)\right)^2 + \frac{1}{2}\exp(-2ct) \sup_{\tau \in [t_0, t]} \left(a_6(2a(2\|u(\tau)\|_U))\right)^2$$

$$+ \frac{1}{2}\exp(-2ct)\left(a_6(4a(4\|x_0\|_X))\right)^2$$

or

$$V(t+1) \leq \exp(-c) V(t) + \exp(-2ct) \sup_{\tau \in [t_0, t]} \rho_1\left(r(\tau)\|u(\tau)\|_U\right) + \exp(-2ct) \rho_2\left(\|x_0\|_X\right) \tag{3.23}$$

where $\rho_1(s) := a_5(s) + a_6(4a(4R))a_5(s) + (a_5(s))^2 + \frac{1}{2}(a_6(2a(2s)))^2$, $\rho_2(s) := \frac{1}{2}(a_6(4a(4s)))^2$ and $r(t) := q(t) + 1$.
Inequality (3.23) in conjunction with (3.18a) directly implies for all $t \geq t_0$:

$$a_1\big(\|H(t,x(t))\|_Y\big) \leq \exp(-c(t-t_0))a_2\big(\widetilde{\beta}(t_0)\|x_0\|_X\big)$$
$$+ \exp(-c(t-t_0)) \sup_{\tau \in [t_0,t]} \frac{\exp(2c)}{\exp(c)-1} \rho_1\big(r(\tau)\|u(\tau)\|_U\big) + \exp(-c(t-t_0)) \frac{\exp(2c)}{\exp(c)-1} \rho_2(\|x_0\|_X) \quad (3.24)$$

Finally, inequality (3.24) implies inequality (3.3) with $\sigma(s,t) := 2a_1^{-1}\left(2\exp(-ct)a_2(s) + 2\exp(-ct)\frac{\exp(2c)}{\exp(c)-1}\rho_2(s)\right)$,
$\beta(t) := \widetilde{\beta}(t) + 1$, $\delta(t) := r(t)$ and $\zeta(s) := 2\frac{\exp(2c)}{\exp(c)-1}\rho_1(s)$. The proof is complete. ◁

**Example 3.4** Consider the nonlinear finite-dimensional discrete-time time-varying system:

$$\begin{aligned}
x_1(t+1) &= d(t)x_1(t) \\
x_2(t+1) &= 2^{-t}d(t)|x_1(t)|^{\frac{1}{2}} + u(t) \\
Y(t) &= H(t,x(t)) := x_2(t) \\
x(t) &:= (x_1(t), x_2(t)) \in \Re^2, \, t \in Z^+, \, d(t) \in [-2,2], \, u(t) \in \Re
\end{aligned} \quad (3.25)$$

Since the corresponding "unforced" system (3.25) with $u(t) \equiv 0$ coincides with system (2.18), which was studied in Example 2.3 and was proved to be non-uniformly in time RGAOS, we conclude by virtue of Proposition 3.3 that system (3.25) satisfies the non-uniform in time IOS property. In order to determine the functions $\sigma \in KL$, $\beta, \gamma \in K^+$ and $\rho \in K_\infty$ for which (3.2) is satisfied, we have to consider the continuous function $V(t,x) := \exp(-t)|x_1| + |x_2|$ (defined in Example 2.3), which clearly satisfies the following inequality for all $(t,x,d,u) \in Z^+ \times \Re^2 \times [-2,2] \times \Re$:

$$V(t+1, dx_1, 2^{-t}d|x_1|^{\frac{1}{2}} + u) \leq \exp(-t-1)|d||x_1| + 2^{-t}|d||x_1|^{\frac{1}{2}} + |u| \leq 2e^{-1}\exp(-t)|x_1| + 2^{-t+1}|x_1|^{\frac{1}{2}} + |u| \quad (3.26)$$

Let $x(t)$ denote the unique solution of (3.25) initiated from $x_0 \in \Re^2$ at time $t_0 \in Z^+$ and corresponding to $(d,u) \in M_{[-2,2]} \times M_\Re$. It can be easily shown (using induction) that the component $x_1(t)$ of the solution satisfies the estimate $|x_1(t)| \leq 2^{(t-t_0)}|x_0|$ for all $t \geq t_0$. Let $V(t) := V(t, x(t))$ and notice that inequality (3.26) in conjunction with the previous estimate for $x_1(t)$ gives:

$$V(t+1) \leq 2e^{-1}V(t) + 2^{-\frac{t}{2}+1}|x_0|^{\frac{1}{2}} + |u(t)|, \quad \forall t \geq t_0 \quad (3.27)$$

Using induction arguments and inequality (3.27), we obtain the following estimate for $V(t)$:

$$V(t) \leq \exp(-c(t-t_0))\left(V(t_0) + \frac{4}{1-2e^{-1}}|x_0|\right) + \frac{2}{1-2e^{-1}} \sup_{t_0 \leq \tau \leq t}(\exp(-c(t-\tau))|u(\tau)|)$$

where $c := \log\left(\frac{2}{1+2e^{-1}}\right)$. The latter inequality combined with (2.19) implies (3.2) with $\sigma(s,t) := 6K\exp(-ct)s$, $\beta(t) \equiv 1$, $\rho(s) := \frac{1}{3}s$ and $\gamma(t) \equiv 1$, where $K := \frac{1}{1-2e^{-1}}$. ◁

## 4. The Robust Output Feedback Stabilization (ROFS) Problem

In this section we first introduce the reader to the notion of the ROFS problem for discrete-time systems (see [13] for the ROFS problem for continuous-time systems). Throughout this section we make the following technical assumption for the "measured output" map $h: Z^+ \times \mathcal{X} \to Y'$ involved in (1.4):

**(A3)** The output map $h \in CU(Z^+ \times \mathcal{X}; Y')$ involved in (1.4), with $h(t,0) = 0$ for all $t \in Z^+$, satisfies:
  (1) There exists a set $S \subseteq Y'$ such that $S = h(t, \mathcal{X})$ for all $t \in Z^+$.
  (2) There exists a function $a \in CU(Y'; S)$, such that for every $y \in S$, it holds that $a(y) = y$.

**Definition 4.1** *Consider system (1.2a) with output maps given by (1.2b) and (1.4) under hypotheses (A2-3) and (H2). The output $Y = H(t,x)$ is called the "stabilized output" while the output $y = h(t,x)$ is called the "measured output".*

**1)** *The problem of **continuous static Robust Output Feedback Stabilization (continuous static ROFS)** for (1.2) with measured output $y = h(t,x)$ and stabilized output $Y = H(t,x)$ is said to be **globally solvable** if there exists a **continuous** function $k \in CU(Z^+ \times S; U)$ (where $S$ is the set involved in hypothesis (A3)) with $f(t,d,0,k(t,0)) = 0$ for all $(t,d) \in Z^+ \times D$, such that the closed-loop system (1.2a,b) with $u(t) = k(t, h(t, x(t)))$ is non-uniformly in time RGAOS. Particularly, we say that the feedback function $k \in CU(Z^+ \times S; U)$ **globally solves** the continuous static ROFS problem for (1.2) with measured output $y = h(t,x)$ and stabilized output $Y = H(t,x)$. Moreover, if the set $H^{-1}(t,0)$ is positively invariant for the closed-loop system (1.2a,b) with $u(t) = k(t, h(t, x(t)))$, then we say that the continuous static ROFS problem for (1.2) with measured output $y = h(t,x)$ and stabilized output $Y = H(t,x)$ is **globally strongly solvable**.*

**2)** *The problem of **continuous dynamic Robust Output Feedback Stabilization (continuous dynamic ROFS)** for (1.2) with measured output $y = h(t,x)$ and stabilized output $Y = H(t,x)$ is said to be **globally solvable** if there exist a normed linear space $W$, **continuous** functions $k \in CU(Z^+ \times S \times W; U)$, $g \in CU(Z^+ \times S \times W; W)$ with $f(t,d,0,k(t,0,0)) = 0$, $g(t,0,0) = 0$ for all $(t,d) \in Z^+ \times D$, such that the following system with state space $\mathcal{X} \times W$ is non-uniformly in time RGAOS:*

$$\begin{aligned} x(t+1) &= f(t, d(t), x(t), k(t, h(t, x(t)), w(t))) \\ w(t+1) &= g(t, h(t, x(t)), w(t)) \\ Y(t) &= H(t, x(t)) \end{aligned} \quad (4.1)$$

*Moreover, if the set $H^{-1}(t,0) \times W$ is positively invariant for system (4.1) then we say the continuous dynamic ROFS problem for (1.2) with measured output $y = h(t,x)$ and stabilized output $Y = H(t,x)$ is **globally strongly solvable.***

The following result is an immediate consequence of Proposition 2.2 and provides necessary Lyapunov-like conditions for the solvability of the static ROFS problem.

**Theorem 4.2** *Consider the ROFS problem for (1.2) with measured output given by (1.4) under hypotheses (A2-3) and (H2). Suppose that the continuous static ROFS problem for (1.2) with measured output $y = h(t,x)$ and stabilized output $Y = H(t,x)$ is globally solvable. Then there exist functions $V \in CU(Z^+ \times \mathcal{X}; \Re^+)$, $a_1, a_2 \in K_\infty$, $\beta_1, \beta_2 \in K^+$ and constant $\lambda \in (0,1)$ such that the following inequalities hold:*

$$a_1\left(\|H(t,x)\|_Y\right) \leq V(t,x) \leq a_2\left(\beta_2(t)\|x\|_\mathcal{X}\right), \ \forall (t,x) \in Z^+ \times \mathcal{X} \quad (4.2a)$$

$$\inf_{u \in U} \sup\left\{V(t+1, f(t,d,x,u)) - \lambda V(t,x); x \in h^{-1}(t,y), d \in D\right\} \leq 0, \ \forall (t,y) \in Z^+ \times S \quad (4.2b)$$

*If the static ROFS problem for (1.2) with measured output $y = h(t,x)$ and stabilized output $Y = H(t,x)$ is globally strongly solvable, then the following condition is additionally satisfied:*

$$\inf_{u \in U^*(t,y)} \sup\left\{V(t+1, f(t,d,x,u)) - \lambda V(t,x); x \in h^{-1}(t,y), d \in D\right\} \leq 0$$

$$\forall t \geq 0 \text{ and for all } y \in S \text{ for which } H^{-1}(t,0) \cap h^{-1}(t,y) \neq \varnothing \quad (4.2c)$$

where
$$U^*(t,y) := \{u \in U \,;\, H(t+1, f(t,d,x,u)) = 0 \text{ for all } d \in D \text{ and } x \in H^{-1}(t,0) \cap h^{-1}(t,y)\} \quad (4.2d)$$

*Finally, if the feedback function* $k \in CU(Z^+ \times S; U)$ *that globally solves the static ROFS problem for (1.2) with measured output* $y = h(t,x)$ *and stabilized output* $Y = H(t,x)$ *satisfies* $k(t,0) = 0$ *for all* $t \in Z^+$, *then the following condition is satisfied:*
$$\sup\{V(t+1, f(t,d,x,0)) - \lambda V(t,x) \,;\, x \in h^{-1}(t,0), d \in D\} \leq 0, \quad \forall t \in Z^+ \quad (4.2e)$$

We also notice the following fact, which combined with Theorem 4.2 provides necessary conditions for the solvability of the continuous dynamic ROFS problem:

**Fact V:** *The continuous dynamic ROFS problem for (1.2) with measured output* $y = h(t,x)$ *and stabilized output* $Y = H(t,x)$ *is globally (strongly) solvable if and only if the continuous static ROFS problem for the following system:*

$$\begin{aligned} x(t+1) &= f(t, d(t), x(t), u(t)) \\ w(t+1) &= v(t) \\ (x(t), w(t)) &\in X \times W, \, t \in Z^+, \, d(t) \in D, \, (u(t), v(t)) \in U \times W \end{aligned} \quad (4.3)$$

*with stabilized output* $Y = H(t,x)$ *and measured output* $\tilde{y} = (h(t,x), w)$ *is globally (strongly) solvable.*

We next give the notion of robust complete observability for discrete-time systems. The definition given here directly extends the corresponding notions given in [3,26], concerning autonomous continuous-time systems.

**Definition 4.3** *Consider the system (1.2a) and let* $(d_i, u_i) \in D \times U$, $i = 0,1,...$ *and define recursively the following family of continuous mappings:*
$$F_0(t,x) = x, \, F_1(t,x,d^{(1)},u^{(1)}) = f(t,d_0,x,u_0)$$
$$F_i(t,x,d^{(i)},u^{(i)}) := f(t+i-1, d_{i-1}, F_{i-1}(t,x,d^{(i-1)},u^{(i-1)}), u_{i-1}), \, i \geq 2$$
$$y_0(t,x) = h(t,x), \, y_i(t,x,d^{(i)},u^{(i)}) := h(t+i, F_i(t,x,d^{(i)},u^{(i)})), \, i \geq 1$$

*where* $d^{(i)} := (d_0,...,d_{i-1})$, $u^{(i)} := (u_0,...,u_{i-1})$ *for* $i \geq 1$. *Let an integer* $p \geq 1$ *and define the continuous mapping for all* $(t,x,d^{(p)},u^{(p)}) \in \Re^+ \times X \times D^p \times U^p$:

$$y^{(p)}(t,x,d^{(p)},u^{(p)}) := (y_0(t,x),..., y_{p-1}(t,x,d^{(p-1)},u^{(p-1)}))$$

*We say that a continuous function* $k \in CU(Z^+ \times X; W)$, *where* $W$ *is a normed linear space, is **robustly completely observable from the output*** $y = h(t,x)$ ***with respect to (1.2a)*** *if there exists an integer* $p \geq 1$ *and a continuous function (called the reconstruction map)* $\Psi \in CU(Z^+ \times S \times S^p \times U^p; W)$ *such that for all* $(t,x,d^{(p)},u^{(p)}) \in Z^+ \times X \times D^p \times U^p$ *it holds that*

$$k(t+p, F_p(t,x,d^{(p)},u^{(p)})) = \Psi(t+p, y_p(t,x,d^{(p)},u^{(p)}), y^{(p)}(t,x,d^{(p)},u^{(p)}), u^{(p)}) \quad (4.4)$$

*We say that system (1.2a) is **robustly completely observable from the output*** $y = h(t,x)$ *if the identity function* $k(t,x) = x$ *is completely observable.*

**Remark 4.4: (a)** Notice that for every input $(d,u) \in M_D \times M_U$ and for every $(t_0, x_0) \in Z^+ \times X$, the unique solution $x(t)$ of (1.2a) corresponding to $(d,u)$ and initiated from $x_0$ at time $t_0$, satisfies the following relation:

$$k(t, x(t)) = \Psi(t, y(t), y(t-p), y(t-p+1),..., y(t-1), u(t-p),..., u(t-1)), \, \forall t \geq t_0 + p$$

Following the terminology in [22], if system (1.2a) is robustly completely observable from the output $y = h(t,x)$ then every control $(d,u) \in M_D \times M_U$ final-state distinguishes between any two events in time $p \in Z^+$.

**(b)** Notice that every continuous function of the measured output $k(t,x) = \theta(t, h(t,x))$, where $\theta : Z^+ \times S \to W$ is a continuous function with the following property:

"for every pair of bounded sets $I \subset Z^+$, $A \subseteq S$ and for every $\varepsilon > 0$ the set $\theta(I \times A)$ is bounded and there exists $\delta > 0$ such that $\|\theta(t,y) - \theta(t,y_0)\|_W < \varepsilon$, for all $t \in I$, $y, y_0 \in A$ with $\|y - y_0\|_{Y'} < \delta$" is robustly completely observable from the measured output.

**(c)** Notice that since $0 \in X$ is an equilibrium point for (1.2) and $h(t,0) = 0$ for all $t \geq 0$, by setting $x = 0$, $u^{(p)} = 0$ in (4.4), we obtain:
$$k(t,0) = \Psi(t,0,0), \quad \forall t \geq p$$

Without loss of generality we may assume that the reconstruction map $\Psi$ is continuously extended to $Z^+ \times S \times S^p \times U^p$ so that the above equality holds for all $t \in Z^+$.

The following proposition provides sufficient conditions for the solvability of the ROFS problem for (1.2).

**Proposition 4.5** *Consider the ROFS problem for (1.2) with measured output given by (1.4) under hypotheses (A2-3) and (H2). Suppose that:*

*(i) There exists a continuous function $k \in CU(Z^+ \times X; U)$ with $f(t,d,0,k(t,0)) = 0$ for all $(t,d) \in Z^+ \times D$, such that the closed-loop system (1.2a,b) with $u(t) = k(t,x(t))$ is non-uniformly in time RGAOS.*

*(ii) The feedback function $k \in CU(Z^+ \times X; U)$ is robustly completely observable from the output $y = h(t,x)$ with respect to (1.2a).*

*Then the continuous dynamic ROFS problem for (1.2) with measured output $y = h(t,x)$ and stabilized output $Y = H(t,x)$ is globally solvable.*

**Proof** Since $k \in CU(Z^+ \times X; U)$ is robustly completely observable from the output $y = h(t,x)$ with respect to (1.2a) there exists an integer $p \geq 1$ and a reconstruction map $\Psi \in CU(\Re^+ \times S \times S^p \times U^p; U)$ such that for all $(t,x,d^{(p)},u^{(p)}) \in Z^+ \times X \times D^p \times U^p$ (4.4) holds. Consider the following system:

$$\begin{aligned}
w_1(t+1) &= y(t) \\
w_2(t+1) &= w_1(t) \\
&\vdots \\
w_p(t+1) &= w_{p-1}(t) \\
w_{p+1}(t+1) &= u(t) \\
w_{p+2}(t+1) &= w_{p+1}(t) \\
&\vdots \\
w_{2p}(t+1) &= w_{2p-1}(t) \\
u(t) &= \Psi(t, y(t), P_S(w(t))) \\
w_i(t) &\in S \subseteq Y', \, i = 1,\ldots,p \\
w_i(t) &\in U, \, i = p+1,\ldots,2p \\
w(t) &:= (w_1(t),\ldots,w_{2p}(t)) \in W := (Y')^p \times U^p, \, t \in Z^+
\end{aligned} \quad (4.5)$$

where
$$P_S(w) := (a(w_p), a(w_{p-1}), \ldots, a(w_1), w_{2p}, w_{2p-2}, \ldots, w_{p+1}) \quad (4.6)$$

and $a: Y' \to S$ is the continuous function involved in hypothesis (A3). Clearly, for every $(t_0, x_0, w_0, d) \in Z^+ \times X \times W \times M_D$ the solution of (1.2) with (4.5) and initial condition $(x(t_0), w(t_0)) = (x_0, w_0)$ corresponding to input $d \in M_D$ satisfies for all $t \geq t_0 + p$:

$$w_i(t) = y(t-i) \quad , \quad i = 1, \ldots, p$$
$$w_{p+i}(t) = u(t-i) \quad , \quad i = 1, \ldots, p \tag{4.7}$$

and consequently by virtue of hypothesis (A3) and definition (4.6) we obtain that:

$$P_S(w(t)) = (w_p(t), w_{p-1}(t), \ldots, w_1(t), w_{2p}(t), w_{2p-1}(t), \ldots, w_{p+1}(t)), \text{ for all } t \geq t_0 + p \tag{4.8}$$

It follows from (4.7) and (4.8) and Remark 4.4 that:

$$u(t) = k(t, x(t)) = \Psi(t, y(t), P_S(w(t))), \quad \forall t \geq t_0 + p \tag{4.9}$$

Equality (4.9) shows that the implemented control action given by $u(t) = \Psi(t, y(t), P_S(w(t)))$ coincides with the control action given by the state feedback law $u(t) = k(t, x(t))$ after $p$ time units. Since the closed-loop system (1.2a,b) with $u(t) = k(t, x(t))$ is non-uniformly in time RGAOS, there exist functions $\sigma \in KL$, $\beta \in K^+$ such that for every $d \in M_D$, $(t_0, x_0) \in Z^+ \times X$, the unique solution $x(t)$ of (1.2a,b) with $u(t) = k(t, x(t))$ initiated from $x_0 \in X$ at time $t_0 \in Z^+$ and corresponding to $d \in M_D$, satisfies (2.2). It follows from (4.9) that for every $d \in M_D$, $(t_0, x_0, w_0) \in Z^+ \times X \times W$, the unique solution $(x(t), w(t))$ of (1.2) with (4.5) and initial condition $(x(t_0), w(t_0)) = (x_0, w_0)$ corresponding to input $d \in M_D$ satisfies

$$\|H(t, x(t))\|_Y \leq \sigma\big(\beta(t_0 + p)\|x(t_0 + p)\|_X, t - t_0 - p\big), \quad \forall t \geq t_0 + p \tag{4.10}$$

Notice that by virtue of Remark 4.4 (c) and since $f(t, d, 0, k(t, 0)) = 0$ for all $(t, d) \in Z^+ \times D$, we may conclude that $0 \in X \times W$ is an equilibrium point for system (1.2) with (4.5). Moreover, by virtue of hypotheses (A2-3) and (H2), it follows that system (1.2) with (4.5) satisfies hypotheses (H2-3) and consequently by virtue of Fact I, there exist functions $\mu \in K^+$, $a \in K_\infty$, such that for every $d \in M_D$, $(t_0, x_0, w_0) \in Z^+ \times X \times W$, the unique solution $x(t)$ of (1.2) with (4.5) initiated from $(x(t_0), w(t_0)) = (x_0, w_0)$ at time $t_0 \geq 0$ and corresponding to $d \in M_D$, satisfies:

$$\|x(t)\|_X + \|w(t)\|_W \leq \mu(t) \, a\big(\|x_0\|_X + \|w_0\|_W\big), \quad \forall t \geq t_0 \tag{4.11}$$

Combining estimates (4.10) and (4.11) we conclude that the closed-loop system (1.2) with (4.5) satisfies the Robust Output Attractivity Property (property P2 of Definition 2.1). By virtue of Fact II, the closed-loop system (1.2) with (4.5) is non-uniformly in time RGAOS. The proof is complete.   ◁

An immediate consequence of Proposition 4.5 is the following proposition, which provides a necessary and sufficient condition for the solvability of the dynamic ROFS problem for (1.2).

**Proposition 4.6 (Separation Principle)** *Consider the ROFS problem for (1.2) with measured output given by (1.4) under hypotheses (A2-3) and (H2). The following statements are equivalent:*

**(a)** *There exist a normed linear space $W'$, continuous functions $k \in CU(Z^+ \times X \times W'; U)$, $g \in CU(Z^+ \times X \times W'; W')$ with $f(t, d, 0, k(t, 0, 0)) = 0$, $g(t, 0, 0) = 0$ for all $(t, d) \in Z^+ \times D$, such that the following system with state space $X \times W'$ is non-uniformly in time RGAOS:*

$$\begin{aligned} x(t+1) &= f(t, d(t), x(t), k(t, x(t), w'(t))) \\ w'(t+1) &= g(t, x(t), w'(t)) \\ Y(t) &= H(t, x(t)) \end{aligned} \tag{4.12}$$

*Moreover, the functions $k \in CU(Z^+ \times X \times W'; U)$ and $g \in CU(Z^+ \times X \times W'; W')$ are robustly completely observable from the output $y' = (h(t, x), w')$ with respect to the system:*

$$x(t+1) = f(t, d(t), x(t), u_1(t))$$
$$w'(t+1) = u_2(t) \qquad (4.13)$$
$$(x(t), w'(t)) \in X \times W', \, u(t) = (u_1(t), u_2(t)) \in U \times W', \, d(t) \in D, \, t \in Z^+$$

**(b)** *The continuous dynamic ROFS problem for (1.2) with measured output $y = h(t, x)$ and stabilized output $Y = H(t, x)$ is globally solvable.*

Implication (a) $\Rightarrow$ (b) of Proposition 4.6 is an immediate application of Proposition 4.5 to the control system (4.13) with input $(u_1, u_2)$. Implication (b) $\Rightarrow$ (a) of Proposition 4.6 is an immediate consequence of Remark 4.4(b) and Definitions 4.1, 4.3. We remark that since the component $x(t)$ of the solution of (4.13) does not depend on the input $u_2$, the requirement that the functions $k \in CU(Z^+ \times X \times W'; U)$ and $g \in CU(Z^+ \times X \times W'; W')$ are robustly completely observable from the output $y' = (h(t, x), w')$ with respect to the system (4.13) implies the requirement that for every $w' \in W'$ the functions $(t, x) \in Z^+ \times X \to k(t, x, w') \in U$ and $(t, x) \in Z^+ \times X \to g(t, x, w') \in W'$ are robustly completely observable from the output $y = h(t, x)$ with respect to the system (1.2).

**Example 4.7** Consider the ROFS problem for the system

$$\begin{aligned} x_1(t+1) &= x_2(t) \\ x_2(t+1) &= x_2^2(t) + u(t) \\ x_3(t+1) &= d(t)x_3(t) + \exp(t)x_2(t) \\ Y(t) &= x(t) \\ x &:= (x_1, x_2, x_3) \in \Re^3, \, u(t) \in \Re, \, t \in Z^+, \, d(t) \in [-r, r] \end{aligned} \qquad (4.14)$$

where $r \in [0, 1)$, with measured output $y = x_1$. First notice that the feedback function $k(t, x) := -x_2^2$ stabilizes system (4.14), non-uniformly in time. We prove this claim by considering the Lyapunov function $V(t, x) := |x_1| + 3 \exp(t)|x_2| + |x_3|$, which clearly satisfies the following inequalities:

$$|Y| = |x| \le V(t, x) \le 5 \exp(t)|x|, \, \forall (t, x) \in Z^+ \times \Re^3 \qquad (4.15a)$$

$$V(t+1, x_2, x_2^2 + k(t, x), dx_3 + \exp(t)x_2) \le (1 + \exp(t))|x_2| + r|x_3| \le \max\left\{\frac{2}{3}, r\right\} V(t, x), \, \forall (t, x, d) \in Z^+ \times \Re^3 \times [-r, r] \qquad (4.15b)$$

and since $\max\left\{\frac{2}{3}, r\right\} < 1$, by virtue of Proposition 2.2, we conclude that the closed-loop system (4.14) with $u(t) = k(t, x(t))$ is non-uniformly in time RGAOS. Moreover, the feedback function $k(t, x) := -x_2^2$ is robustly completely observable from the output $y = x_1$. Particularly, we define the continuous mappings (following the notation of Definition 4.3):

$$F_0(t, x) = x, \, F_1(t, x, d^{(1)}, u^{(1)}) = \begin{pmatrix} x_2 \\ x_2^2 + u_0 \\ d_0 x_3 + \exp(t) x_2 \end{pmatrix}, \, y_0(t, x) = x_1, \, y_1(t, x, d^{(1)}, u^{(1)}) := x_2$$

Clearly, we have:

$$k\left(t+1, F_1(t, x, d^{(1)}, u^{(1)})\right) = \Psi(t+1, y_1, y_0, u_0) := -\left(y_1^2 + u_0\right)^2$$

Consequently, the closed-loop system (4.14) with:

$$\begin{aligned} w(t+1) &= u(t) \\ u(t) &= -\left(y^2(t) + w(t)\right)^2 \\ w(t) &\in \Re, \, t \in Z^+ \end{aligned}$$

is non-uniformly in time RGAOS. ◁

## 5. Conclusions

The notions of non-uniform in time Robust Global Asymptotic Output Stability (RGAOS) and non-uniform in time Input-to-Output Stability (IOS) are studied in the present paper for time-varying discrete-time systems. Characterizations and links between these notions are provided. Particularly, it is shown that a discrete-time system with continuous dynamics satisfies the non-uniform in time IOS property if and only if the corresponding unforced system is non-uniformly in time RGAOS. The Robust Output Feedback Stabilization (ROFS) problem is studied next. Necessary and sufficient conditions for the solvability of this problem are provided.